\input amstex
\documentstyle{amsppt}
\magnification=1200 
\parindent 20 pt
\NoBlackBoxes

\define\iy{\infty}

\define\a{\alpha}

\define\3{\char\ss}
\define\noin{\noindent}
\define\r{\rightarrow}

\define\s{\subset}

\define\BC{\Bbb C}
\define\BR{\Bbb R}
\define\BP{\Bbb P}

\define \1{^{-1}}

\define \sbt{\subset}

\define \ri{\rightarrow}
\define \bk{\bigskip} 
\define \mk{\medskip} 
\define \sk{\smallskip}
\define \ep{\endproclaim}
\define \edm{\enddemo}

\define \fc{\frac}

\define\BX{\Bbb X}
 \define\CX{\Cal X}
\define\CZ{\Cal Z}
\define\CY{\Cal Y}

\def\cupl{\operatornamewithlimits{\bigcup}\limits}

\baselineskip 20pt
\title
{ Estimates 
of the number of rational mappings 
from a fixed variety to varieties of general type}
\endtitle
\leftheadtext{T. Bandman and G. Dethloff}
\rightheadtext{Estimates on the number of maps}
\author
T. Bandman and G. Dethloff
\endauthor
\thanks The first named author  was supported by the 
 Ministry of Science and
Technology, by the Ministry of Absorption, State of Israel,
the Emmy N\"other Institute for
Mathematics, Bar-Ilan University and the SFB 170 `Geometry and Analysis',
 G\"ottingen.
The second named author was supported by the MSRI Berkeley and by the DFG.
We would like to thank the referee for his suggestions how to improve the style 
of this paper.
\endthanks

\abstract
First we find  effective bounds  
 for the number of  dominant rational maps $f:X \rightarrow Y$
 between two fixed 
smooth projective varieties with ample canonical bundles. The bounds are
of the type $\{A \cdot K_X^n\}^{\{B \cdot K_X^n\}^2}$, where $n=dimX$, $K_X$ is 
the canonical bundle of $X$ and $A,B $ are some constants, depending only on 
$n$.

Then we show that for any variety $X$   
there exist numbers $c(X)$ and $C(X)$ with the following
 properties:

 For any
threefold $Y$ of general type the number of dominant rational maps
$f:X \r Y$ is bounded above by $c(X)$.

 The number
of threefolds $Y$, modulo birational equivalence,
 for which there exist dominant rational maps $f:X \r Y$, is bounded above by
$C(X)$.
 
\noindent If, moreover, $X$ is a threefold of general type, we prove that
$c(X)$ and $C(X)$ only depend on the index $r_{X_c}$ of the canonical model
$X_c$ of $X$ and on $K_{X_c}^3$.
\endabstract
\endtopmatter 
\document
\mk
\subheading{0. Introduction}

\sk

 Let $X$ and $Y$ be algebraic varieties, i.e. complete integral schemes
over a field of characteristic zero.  Denote by $R(X,Y)$ the set of  
dominant rational maps $f:X \r Y$.
Then the classical theorems of de Franchis \cite{Fra} and Severi 
(cf. \cite{Sam}) can be stated as
follows:

\proclaim{Theorem 0.1 }\

a) (de Franchis): For any Riemann surface $X$ and any hyperbolic 
Riemann surface $Y$
the set $R(X,Y)$ is finite. Furthermore, there exists an upper bound
for $\#R(X,Y)$ only in terms of $X$.

b) (Severi) For a fixed algebraic variety $X$ there exist only finitely 
many hyperbolic Riemann surfaces $Y$
such that $R(X,Y)$ is nonempty. 
\ep

 S.
Kobayashi and T. Ochiai  \cite{Kob-Och}  prove the following 
generalization of the de Franchis Theorem: If $X$ is a 
Moishezon space and $Y$
a compact complex space of general type, then the set of surjective
 meromorphic
maps from $X$ to $Y$ is finite.
Other generalizations can be found in \cite{Des-Men1}, \cite {Nog},
\cite {Suz}, and in the survey \cite {Zai-Lin}. 
Generalizations of the second part of de Franchis' Theorem
were given in 
 \cite{Ban1}, \cite{Ban2} and \cite{Ban-Mar}.
In the latter paper it is proved  that
for complex  projective varieties $X$ and $Y$ with only canonical singularities
and nef and big canonical classes $K_X$ and $K_Y$ respectively 
the number
$\#R(X,Y)$ can be bounded in terms of the selfintersection $K_X^3$
of the canonical class of $X$ and of the indices of $X$ and 
$Y$.
   This bound  is not effective. Effective bounds are known only if the 
 varieties $Y$ are curves or surfaces (\cite {Kan}, \cite {How-Som2},
\cite {Tsa3}).

\bigskip
Section 1 of this paper     
   contains an effective estimate of 
   the number of  mappings  in $R(X,Y)$, provided that both varieties 
   $X$ and $Y$ are smooth projective with ample canonical bundles
   $K_X, K_Y$ respectively. This bound has the form
$\{A \cdot K_X^n\}^{\{B \cdot K_X^n\}^2}$, where $n=dimX$, $K_X$ is 
the canonical bundle of $X$ and $A,B $ are some constants, depending only on 
$n$.

 This bound seems to be very big. 
But
 it is known that the bound cannot be
polynomial in $K^n_X$ (\cite {Kan}).
 Moreover, even for the case of curves of genus $5$
the best bound in \cite {Kan} is of order $exp (30).$

The idea to obtain this bound was used in
 (\cite {How-Som1}) for proving
 finiteness of the automorphism
 group of a projective variety with ample
canonical bundle. It could not be made effective at that time, as no effective
variants of the Big Matsusaka theorem were available. 
Moreover,  exponential bounds for the number of automorphisms are not
interesting, as  they should
 be  linear in $K^n_X$
(\cite {Sza}).

\bigskip

 In  sections 2, 3 and 4  we generalize
 Severi's result (Theorem 0.1(b)) to higher dimensions.

Denote by $\Cal {F}(X)$ the set of  pairs $(Y,f)$, where $Y$ is of general
type and $f \in R(X,Y)$. Let $\Cal {F}_m(X) \s \Cal {F}(X) $
the subset of those pairs $(Y,f)$  for which the $m$-th pluricanonical
mapping of a desingularization of $Y$ is birational onto its image.
 Consider the equivalence relations on
 $\Cal {F}$ and $\Cal {F}_m$: $(f:X \r Y) \sim (f_1:X \r Y_1)$
iff  $b \circ f = f_1$, where $b \in $Bir$(Y,Y_1)$.
    The elements of  $\Cal {F}(X)/ \sim$ we call {\it targets}.

The following conjecture is stated by Maehara (\cite{Mae3})  
as Iitaka's Conjecture based on Severi's Theorem:

\sk
\proclaim{Conjecture 0.2}
The set $\Cal {F}(X)/ \sim$ of targets is a finite set.
\ep
\sk

\noindent
Maehara proved in Proposition 6.5 in \cite{Mae2} that in characteristic zero 
 $\Cal {F}_m(X) / \sim$
is finite for all $m$. In 
particular
the Conjecture is valid for surfaces $Y$ (take $m=5$).
Special cases and related aspects are discussed in
[Tsa1] - [Tsa-3], [Des-Men2], [Des-Men3], [Mae1] and [How-Som2].
 
\sk
In section 3 we prove Conjecture 0.2 for the case that the targets are 
complex threefolds (Theorem 3.1).

For the proof we use the following theorem of
Luo \cite{Luo1}, \cite {Det}:

\sk
\proclaim{Theorem 0.3 }\  Consider the set of smooth threefolds $Y$ of general 
type, and denote by $ \chi (Y,\Cal {O}_Y)$ the holomorphic Euler characteristic.
Then for any fixed $\chi = \chi (Y,\Cal {O}_Y)$, there is
a universal integer $m'$ such that $h^0(Y,\Cal {O}_Y(m'K_Y)) \geq 2$.
 Furthermore, there is a universal integer
$m$  such that the $m$-th pluricanonical map 
$\Phi_{mK}: Y \r \Phi_{mK}(Y)$  maps birationally onto its image.
\ep
\sk

In section 4 the domain is a 
threefold of general type. In this case
we show (Theorem 4.1) that there is  a  bound for the number of targets $\# 
\Cal {F}(X)/ \sim$, which depends 
only on the selfintersection $K_{X_c}^3$
and the index $r_{X_c}$ of the canonical model $X_c$  of $X$.

The proof 
is based on the fact, due to Kollar  \cite{Kol1},
that canonical threefolds with fixed Hilbert polynomial form a bounded
family. Using semicontinuity theorems
for the dimensions of cohomology groups we get estimates of the holomorphic
Euler characteristics of the targets.
 Then we show that
 the graphs of maps under consideration, which map from canonical threefolds $X$ with fixed index
$r_X$ and fixed $K_X^3$, form a finite number of algebraic families.
The number of targets  is bounded by the number of 
irreducible components of the members of these families.

\bigskip 
In Section 5 we return to generalizations of de Franchis' result (Theorem 
0.1(a)). Consider a threefold $X$ of general type.
We prove (Theorem 5.1) that there exists a bound
for $\# R(X,Y)$, depending only on $X$.
 Namely, it depends on the selfintersection
$K_{X_c}^3$ and on the index $r_{X_c}$ of the canonical model 
$X_c$ of $X$.

 In the review of  Sh. Kobayashi (\cite {Kob}, problem D3)
the question is raised 
if for a compact complex space $X$ and a hyperbolic compact complex space
$Y$ the number of surjective meromorphic maps from $X$ to $Y$ can be bounded only
in terms of $X$. Theorem 5.1 is an answer to this question 
for threefolds of general type.
    
\bigskip
 Further on all  the varieties  are complex; we do not
make difference between line bundles, divisor classes and the divisors
themselves, if  no confusion may arise. 
We fix resp. recall the following notations, which are used in the paper: 

$X$, $Y$ -- complex varieties;
  
$R(X, Y)$    --   the set of  rational dominant maps from $X$
to  $Y$;

$\Cal {F}(X)$ -- the set of  pairs $(Y,f)$, where $Y$ is of general
type and $f \in R(X,Y)$; 
 
$\Cal {F}_m(X)$ -- the subset of those pairs $(Y,f)$  for which the $m$-th pluricanonical

\hskip 1.5 cm mapping of a desingularization of $Y$ is birational onto its image;

 $(f:X \r Y) \sim (f_1:X \r Y_1)$
iff  $b \circ f = f_1$, where $b \in $Bir$(Y,Y_1)$;

$K_X$  -- the canonical sheaf of a variety $X$ with at most canonical 
singularities;

$K_X^n$  -- the $n$-times selfintersection of the class $K_X$, where
 $ n=\dim X$;


$c_i(X)$ -- the $i^{th}$ Chern class of the variety $ X$;


$H^i(X,D) = H^i(X,\Cal {O}_X (D))$; 
$h^i(X,D) = \dim H^i (X,D)$;

$\chi (X,D) = \sum^n_{i=1} (-1)^i h^i(X, D)$;

$X_c$ -- the canonical model of a variety $X$ of general type of
 dim$_{\BC}X \leq 3$;

$r_X$ --  the index of a variety $X$ with at most canonical singularities;

$\Phi_{mK_Y}$ -- the $m$-th pluricanonical map 
from a variety $Y$ with at most

\hskip 1.4 cm canonical singularities.    
 
\vskip 1.0 cm

 \subheading {1. Effective estimates of R(X,Y) for smooth
manifolds $X,Y$ with ample canonical bundles}   
\sk

 The main Theorem of this section is Theorem 1.6 below.
 It provides an effective estimate for  $\#R(X,Y)$ if
  $X,Y$ are smooth manifolds with ample canonical divisors.
\sk

 We first recall some notations and facts about duality: 

\noindent a) A subspace $E \s {\Bbb P}^N$ is called linear if it is the projectivization
of a linear subspace $E^a \s {\Bbb C}^{N+1}$.
Let $\rho :{\Bbb C}^{N+1} \r ({\Bbb C}^{N+1})^*$ denote the canonical isomorphism
between ${\Bbb C}^{N+1}$ and the space $({\Bbb C}^{N+1})^*$ of linear functionals 
on it, which is given by the standard hermitian product on ${\Bbb C}^{N+1}$.
We denote by $({\Bbb P}^N)^*$ resp. $E^*$ the projectivizations of
$({\Bbb C}^{N+1})^*$ resp. of  $\rho (E^a) \s ({\Bbb C}^{N+1})^*$, and call
them the conjugate spaces to ${\Bbb P}^N$ resp. $E$. (We don't use the word
`dual' here in order not to have confusion with the notion of a dual variety which
is defined below.)

\noindent b) We call a rational mapping  $L\: {\Bbb P}^N
\ri {\Bbb P}^M$ linear if it is the projectivization of a linear map
 $L^a\: {\Bbb C}^{N+1} \ri {\Bbb C}^{M+1}$. The projectivization of the induced map
$(L^a)^*\: ({\Bbb C}^{M+1})^* \ri ({\Bbb C}^{N+1})^*$ is denoted by $L^*$ and
called the dual map to $L$.

\noindent c) Let $Z$ be an  $n$-dimensional projective variety embedded 
into the projective space ${\Bbb P}^N.$ 
In any non-singular point $z\in Z$ the projective tangent plane $T_z$
is well defined.  
In the conjugate projective space $({{\Bbb P}^N})^*$ we consider the 
set $Z_0^V$ of all points
$y\in ({{\Bbb P}^N})^*$ such that the corresponding hyperplane $H_y\sbt {\Bbb P}^N$
 contains the
tangent plane $T_z$ to some nonsingular point $z\in Z$. 
We define the dual variety $Z^V$ of the variety $Z$ to be the closure of 
$Z_0^V$ in  the Zarisky
topology.
   
 \noindent d) The dual varieties
have the following fundamental properties (\cite {Del-Kat}): 
\roster 
\item"1)" If $Z$ is nonsingular and Kodaira dimension
 $k(Z) > -\iy,$ then $Z^V$ is irreducible and $\text{codim}\ Z^V=1;$ 
\item"2)" Moreover, if $L$ is a hyperplane section of $Z\sbt {\Bbb P}^N$,   
 the degree $\deg Z^V$ of the variety $Z^V$ may
 be computed by the Chern classes:
$$\deg Z^V = \sum^n_{i=0}
 (-1)^{n+i} (1+i) c_1^i (L) c_{n-i} (Z), \tag1$$
where $c_1(L)$ is the first Chern class of the line bundle corresponding to $L$;
\item"3)"  $Z^{VV} = Z.$ \endroster

\sk 

Let $X$ and $Y$ be two smooth projective $n$-dimensional varieties
with Kodaira dimension bigger than infinity, and let $E$ and $F$
be  very ample line bundles on $X$ and $Y$ respectively.
Then the varieties $X$ and $Y$ are canonically embedded into 
the projectivizations  of the conjugate spaces to
$H^0(X,E)$ and to $H^0(Y,F)$ respectively, which we denote by 
 ${\Bbb P}^N$ and
${\Bbb P}^M$.

 Let $f\:
X\dashrightarrow Y$ be a rational dominant mapping,  
and let $\Psi : H^0 (Y,F)\ri H^0 (X,E)$ be an injective linear map.
We call $f$ to be induced by the map
$\Psi$ if the projectivization of the dual map $\Psi^*$, restricted
to $X$, is  the map $f$. 
Denote by $R(X,E,Y,F)$ the set 
 $$R(X,E,Y,F) = \{ f \in R(X,Y): 
f \text{ is induced by an injective linear map }$$
 $$\Psi \: H^0(Y,F)\ri H^0(X,E) \}$$
\sk
\proclaim {Proposition 1.1 }\ If the set $R(X,E,Y,F)$ is finite, 
 we have:
$$\# R(X,E,Y,F) \leq m^{\psi(E)}$$
where $$m=\sum\limits^n_{i=0} (-1)^{n+i}(1+i) c_1^i (E) c_{n-i} (X) ,$$
$$\psi (E) = (h^0 (X,E))^2 - 1.$$\ep
\sk

Before we start with the proof of Proposition 1.1, we need two Lemmas.

\proclaim {Lemma 1.2}  Denote by $G$ the set of all linear injections $A:
(\BP^M)^* \to ({\BP^N})^*$ such that $A(Y^V)\s X^V$. Then $G$ is an 
quasiprojective subset of $\BP ^K, K= (N+1)(M+1)-1$, of degree
$$ \deg G \leq (\deg X^V)^K .\tag2$$\endproclaim

\bk
 {\it Proof of Lemma 1.2.}
  Any element of $G$ is defined by a $(N+1)(M+1)$ matrix $A$, and
its components $(a_{ij})$ may be considered as its coordinates in the
projective space ${\Bbb P}^K,\ K=(N+1)(M+1) - 1.$
Since $\text{codim}\ X^V =1,$ it is defined in $({\Bbb P}^N)^*$ by
 a single equation    
$F(z_0,\dots, z_N)=0$ with $\deg F =\deg X^V.$  
If $ y\in Y^V$ we have $A(y)\in X^V$, and so 
$F (Ay) =0.$
For a fixed point $y$ and a fixed polynomial  $F$ this
 is an equation for the coordinates $a_{ij}$      
in the space ${\Bbb P}^K.$  

 This means that for any finite 
 sequence of points $y_1, ...,y_r, y_i \in Y^V,
$
 the set $G$ is contained in the algebraic set 
$G^{(r)}$, defined by equations

$$\align
&F(Ay_1)=0\\
&\vdots\\
&F(Ay_r) = 0\endalign$$
in the space ${\Bbb P}^K$ .

Choose any point $y_1 \in Y^V$. Suppose that $\overline{G} \neq G^{(1)}$,
where $\overline{G}$ denotes the Zariski closure of $G$ in ${\Bbb P}^K$.
 Then there exists
a point $y_2$ such that for some $A\in G^{(1)}$
$$ F(Ay_2)\neq 0.$$
Define the set $G^{(2)}$ by the pair $y_1, y_2$. It follows that $G^{(2)}\s G^{(1)}$,
and  for some component $C$ of $G^{(1)}$ all components of $G^{(2)}$ which lie in
$C$ (if there are any at all)  are of smaller dimension than $C$.
  After
performing a finite number of such steps we  get a set 
 $y_1, y_2,...,y_r$, such that $G^{(r)}=\overline{G}$. Hence, $G$
can be defined in $\BP^K$ by  equations of  degree
$deg F= deg X^V$ only. Now, the inequality   
$$ \deg G \leq (deg X^V)^K.$$
follows from the 
 \proclaim {Sublemma (the analogue of the Bezout Theorem)}

Let $X\subset \BP^n$ be an irreducible variety, $ \dim X=i, \deg X=a.$ 
Let $F_1,...F_s$ be homogeneous polynomials of degree $d$ and 
$X_s=\{z\in \BP^n: F_1(z)=F_2(z)=...=F_s(z)=0\}$. Assume that 
$X\cap X_s = \cupl_{j=1}^{N}B_j $ is a union of irreducible 
components $B_j.$ Then
$$\deg (X \cap X_s) =\sum_j \deg B_j \leq ad^i.$$
\ep

\demo{Proof of the Sublemma}
We perform induction by $i=\dim X.$ If $i=1,$ 
   there are two possibilities:

   1. $F_k \bigm|_X=0$ for all $k=0,...s;$ then $X=B_1, N=1,$ 
   and $\deg B_1=\deg X=a.$

   2. $F_1\bigm|_X\not\equiv  0.$ Then $X \cap X_s\subset X\cap X_1$
   is a finite number $T$ of points and $$T\leq \deg 
   (X\cap X_1)\leq ad.$$ 

   Assume that the fact is true for every $i<m.$ 
   If $F_k\bigm|_X = 0$ for all $k=1,...,s,$ then 
   $N=1, X=B_1$ and $\deg B_1=a.$
   If $ F_s\bigm|_X\not\equiv  0,$ then
$X\cap\{F_s=0\}=\bigcup A_q$ is a union of 
irreducible components $A_q$ such that $n_q=\dim A_q<m$
and $$\sum_q \deg A_q \leq ad$$
(see, for example, \cite{Har}, Th.7.,ch. 1).
Let $A_q\cap X_{s-1}=\bigcup_r B_q^r.$
Since $$ \bigcup_{q,r} B_q^r =\left(\bigcup_q A_q\right) \cap X_{s-1}=
X\cap X_s= \bigcup_j B_j,$$ and all  
$B_q^r$ and $B_j$ are irreducible, we
obtain that for any $j$ there are numbers (q,r)
such that $B_j=B_q^r.$ Thus $$\sum_j \deg B_j \leq \sum_{q,r} \deg B_q^r.$$
By induction assumption
$$\sum_r \deg B^r_q\leq\deg A_q \  d^{n_q}\leq \deg A_q \ d^{m-1}.$$
 Summation over  $q$ provides the desired inequality:
 $$\sum_j \deg B_j \leq \sum_q \sum_r \deg B_q^r \leq \sum_q\deg A_q d^{m-1}
 \leq d^{m-1}\sum_q\deg A_q\leq ad^m.$$
\qed \edm

Let  $G=\bigcup_iG_i$ be the decomposition of $G$ in irreducible components $G_i$.

\proclaim {Lemma 1.3}
  Suppose that the points $t_1,t_2 \in G_i$ define
linear maps $A_j,\ \ j=1,2$, which
are  dual to  linear projections
$A_j^*:\BP^N \to\BP^M$ satisfying $A_j^*(X)=Y$ (i.e. $f_j :=A_j^*|_X 
\in R(X,E,Y,F)$).
Then $f_1=f_2$.\endproclaim
 \bk

\demo{Proof of Lemma 1.3}
From now on we fix a basis in $\BP^N$ and $\BP^M.$
Let $t\in G_i$ and let $A_t$ be a linear embedding $A_t: (\BP^M)^*\to(\BP^N)^*$
corresponding to a point $t.$

Consider the following diagram:
 $$\alignat7
&X\quad &&\subset\qquad &&\BP^N\qquad &&\simeq\qquad &&(\BP^N)^*\quad&&
\supset\quad
&&X^V\\
\quad\\
&\downarrow \ A_t^*\bigm|_X\quad  && &&\downarrow\ A_t^*&& &&\cup\\
\quad\\
&X_t\quad &&\subset\qquad &&L_t\qquad &&\simeq\qquad &&E_t\quad&& \supset\quad
&&X^V\cap E_t\\
\quad\\
&\underset \vee \to\vdots \ \tau_t \quad  && &&\downarrow \tau_t&&
&&\uparrow\ A_t\\
\quad\\
&Y\quad &&\subset\qquad &&\BP^M\qquad &&\simeq\qquad &&(\BP^M)^*\quad&&
\supset\quad
&&Y^V\\
 \endalignat
$$
In this diagram:
\roster\item"1)" $A_t$ is a linear embedding of $(\BP^M)^*$ into $(\BP^N)^*;$
\item"2)" $E_t=A_t((\BP^M)^*)$ is a linear subspace of $(\BP^N)^*,$ which is
isomorphic to $(\BP^M)^*;$
\item"3)" The dual map $A_t^*: {\Bbb P}^N \r {\Bbb P}^M$ is the composition of
 a projection of $\BP^N$
onto the subspace
$L_t\subset \BP^N,$ which is dual to $E_t$, and of an isomorphism 
$\tau_t : L_t \r {\Bbb P}^M$ (recall that we have chosen a
basis in $\BP^N$
and a dual basis in $(\BP^N)^*)$. The projection we again denote by $A_t^*$.
\item"4)" $X_t=A_t^*(X)\subset L_t.$
\endroster

\noindent Now the linear map $A_t$ induces a dominant rational map
 $f_t: X\to Y,$ iff  $\tau_t (X_t)=Y$.

\sk

 In order to   proceed  with  the proof, the  following two claims are needed.    

\proclaim {Claim 1}
    Let $R_i(x,t)$, $i=1,...,l$,  be a finite number of polynomials in the variable 
$x \in {\Bbb P}^K$ with coefficients which are polynomials in the variable $t \in
T,$ where $T$ is an irreducible projective variety. Let
$V_t=\{ x \in {\Bbb P}^K: R_1(x,t)=...=R_l(x,t)=0\}$. 
Assume that for some 
point $t_0\in T,$ 
$$r(x,t_0)=rank \{ \frac {\partial R_i}{\partial x_j}|_{(x,t_0)}\}=k$$
 for any $x \in V_{t_0}$. Then the set of points $t\in T$, such that
$r<k$ for some $x \in V_t$, is proper and closed in $T$.\endproclaim

{\it Proof.} Consider the sets $A=\{(x,t) \in {\Bbb P}^K\times T: r(x,t)<k \}$ 
and $B=\{ (x,t) \in {\Bbb P}^K\times T : x \in V_t \}$. 
 Since $A \cap B$ is Zariski closed in ${\Bbb P}^K\times T$, its projection
 to $T$ is Zariski closed in $T.$  Since there is at least one point
$t_0$, which does not belong to the image of this projection, it has to be a proper
 Zariski closed subset.\qed

\sk
By the assumption  $G_i$ has two points $t_1,t_2$, 
defining the maps $f_1$,$f_2 \in R(X,E,Y,F)$.
\proclaim {Claim 2}
Let $T' \s G_i$ be the set of points $t\in G_i$, for which 
 the image $X_t$   
of the projection of $X$ into $L_t$ is smooth and of the same dimension as $X$. 
 Then $T'$ is Zarisky open and contains $t_1, t_2$. \endproclaim
{\it Proof.} We apply Claim 1 with $T=G_i$, ${\Bbb P}^K = L_t$ and where
the $R_i(x,t)$, $i=1,...,l$ are the resultants of the polynomial equations 
 of $X$ in $\BP^N$. Then $V_t=X_t$ is smooth in the point $x$ iff the rank 
$r(x,t)$ is maximal. Since $t_1,t_2$ define
the maps $f_1, f_2 \in R(X,E,Y,F)$, the varieties
$X_{t_i}$ are isomorphic to
$Y$ through the maps $\tau_{t_i}$. Especially $X_{t_1}, X_{t_2}$ are smooth
and of the same dimension as $X$. Now Claim 2
follows from Claim 1. 
\qed
\sk

Using Claim 2, we will show that  all $t \in T'$ correspond to maps 
$f_t \in R(X,E,Y,F)$, and moreover, that $f_t$ does not depend on the parameter
$t \in T'$. Especially, we get $f_{t_1}=f_{t_2}$.

Since $X_{t_1}$ is isomorphic to $Y$, we get for the Kodaira dimensions
$k(X_{t_1}) = k(Y) > - \infty$. Applying the invariance of  plurigenera
to the algebraic family with base $T'$ and fiber $X_t$ (\cite {Har}, 9.13, ch.3),
we get that   $k(X_t)> -\infty$ for all $t \in T'$.
  Then $X_t^V$  has to be irreducible hypersurface in $E_t$,
 contained in $X^V\cap E_t$.
Thus $X^V \cap E_t$ contains an  irreducible component
$C_t$ such that $C_t=X_t^V$, and $C_t^V=X_t^{VV}=X_t$.
 Let $B_{t,i}$ be other irreducible components of $X^V\cap E_t.$
Then $B_{t,i}^V \subset X_t$. Since $X_t$ is irreducible, 
we get $\dim B_{t,i}^V < \dim X_t.$  
Thus, the intersection  $X^V\cap E_t$ is a union of irreducible components 
$B_{t,i}$ and  $C_t$, such that 

a)  the 
components  $B_{t,i} \subset E_t$ are dual to some subsets of $X_t$
(actually, the image of singular points of the projection of $X$ to $X_t$)
of  dimension less than $n=\dim X$; 

b) the component $C_t \subset E_t$ is the only
 one which has $n$ - dimensional dual, and $C_t^V=X_t$.

For any point $t\in T'$ the variety $A_t(Y^V)$ is isomorphic to
$Y^V$. Thus $\{A_t(Y^V)\}^V $ is isomorphic to $Y^{VV}=Y$ and, hence, 
it is $n$-dimensional.
On the other hand, $A_t(Y^V)$ is contained in $X^V \cap E_t$ and is
 irreducible and of the same dimension as $X^V \cap E_t$.
 Hence, $A_t(Y^V)=C_t= X_t^V.$

That means that $A_t$ is an isomorphism between $(\BP^M)^*$ and $E_t$ such that
$A_t(Y^V)=X_t^V\subset E_t.$
Then  the dual isomorphism $\tau_t=A_t^*\bigm|_{L_t}: L_t\to \BP^M$ maps
$(X_t)^{VV}= X_t$ onto
$Y^{VV}=Y$ (see \cite{How-Som1}).
It follows that the map $f_t=\tau_t\circ A_t^*\bigm|_X$ belongs to $R(X,E,Y,F)$.
Since the latter set is finite and the family of the projections $A_t^*$,
 which give $f_t:X \r Y$, varies continuously over $T'$, the map $f_t$ does not depend on $t$.
\qed \edm
\sk
\demo{Proof of Proposition 1.1}\  
A mapping $f\in R(X,E,Y,F)$ is,   by definition, induced by a  linear 
injection
   $\Psi\: H^0 (Y,F)\ri H^0(X,E)$. More precisely, if we denote by
$\tilde{f}: {\Bbb P}^N \r {\Bbb P}^M$ the projectivization of the dual
map $\Psi^*$ to $\Psi$, then $f=\tilde{f}|_X$.

Using the surjectivity of the maps $f:X \r Y$ and $\tilde{f}: 
{\Bbb P}^N \r {\Bbb P}^M$, an easy computation yields that
 the dual map $\tilde{f}^*\: ({\BP^M})^*\to ({\BP^N})^*$
maps $Y^V$ into $ X^V$.

By Lemma 1.3 the number 
of mappings in $ R(X,E,Y,F)$  is, at most,   the number
 of irreducible components  in $G$, which obviously does not exceed the
 sum of degrees of the these components.
Since the mappings in $R(X,E,Y,F)$ are induced by  injections of $H^0(Y,F)$
 into $H^0(X,E)$, we have 
$$M =\dim h^0 (Y,F) -1 \leq \dim h^0(X,E) -1  = N.$$ 
By the formula (1) we get, for $\deg X^V$:
$$\deg X^V = \sum^n_{i=0} (-1)^{n+i}(1+i) c_1^i (E)\cdot c_{n-i} (X).$$
               
To obtain the statement of  
Proposition 1.1 it suffices to insert these
 values into  equation  (2).\ \ \ \qed 
\enddemo
\sk

Using  Proposition 1.1 we obtain the following

\sk
\proclaim{Theorem 1.6 }\  Let $X,Y$ be two smooth complex projective varieties with ample
canonical bundles $K_X$ and $ K_Y$.   
Let $R(X,Y)$ be the set of 
dominant rational maps $f\: X\ri Y,$ and let the
divisors $sK_X, sK_Y$ be very ample (for example $s$  may be $2+12n^n$        
(\cite{Dem})). 
Then 
$$\# R(X,Y) \leq
\left[(-1)^n\sum^n_{i=0} (1+i) s^i c_1^i (X)
c_{n-i} (X)\right]^{\left\{
K^n_X\left(\fc{s^n}{n!} -
\fc{s^{n-1}}{2(n-1)!}+q_{n-2}(s)\right) \right
\}^2},$$       where,  for
each $n$, $q_{n-2}(s)$ is a universal polynomial  of degree $n-2$
. 
\ep

\demo{Proof of Theorem 1.6 } 
Let $E= sK_X$, $F= sK_Y$ and $f \in R(X,E,Y,F)$. 
If $f^*$ here  denotes the pull back of pluricanonical forms by the rational
map $f$, we get, by \cite{Iit}, Theorem 5.3, that 
  $f^*\:H^0(Y, msK_Y)\ri H^0(X,msK_X)$
 is
an  injective linear map for any $m\in {\Bbb N}$. (Since the divisors 
$E$ and $F$ are very ample, any rational map 
$f\in R(X,Y)$ is even regular (\cite {Ban1}), but we don't need this fact here.)
 It is easy to see that the map $f$ is
induced by the linear map $f^*$.
That is why
$R(X,Y)=R(X,sK_X,Y,sK_Y)$. Since this set
is finite (\cite {Kob-Och}), we can apply
 Proposition 1.1: 
$$\align
\# R(X,Y) &= \# R(X,sK_X, Y, sK_Y)\\                                  
&\leq \left[\sum^n_{i=0}(-1)^{n+i} (1+i) s^i c_1^i (K_X) c_{n-i}(X)\right]
^{h^0(X,sK_X)^2-1}\endalign$$
  In this  expression we substitute $ -c_1(K_X)$ by
$c_1(X)$ (these numbers are equal). Further, by  the Riemann-Roch
Theorem
and the Vanishing Theorem for ample line bundles
$$h^0(X,sK_X) = \chi(X, sK_X) = K^n_X \left[\fc{s^n}{n!} -
\fc{s^{n-1}}{2(n-1)!}\right] + P(s),$$
where $P(s)=\sum\limits^{n-2}_{i=0} \a_i s^i$ is a  polynomial of degree $n-2$
in $s,$ the coefficients of which are linear combinations of monomials
of the form $c_I(X) = c_{i_1}(X) \dots c_{i_k} (X), i_1 +\cdots +i_k=n.$

According to   (\cite{Ful-Laz}, \cite {Cat-Sch}), there exist  universal 
constants $D_I$, depending only on $n$, such that
$$|c_I(X)| \leq D_IK^n_X.$$
It follows, that there are other universal constants $\tilde {D}_i$, $i=0,...,n-2$,
which depend only on $n$, such that 
$$|P(s)| \leq \sum^{n-2}_{i=0} |\a_i| s^i
 \leq K^n_X \cdot \sum^{n-2}_{i=0}
\tilde {D}_i\cdot s^i.$$
Hence, it is possible to choose
 $$q_{n-2}(s) = \sum^{n-2}_{i=0} \tilde {D}_is^i.$$\ \ \ \qed

\enddemo
\bigskip

\subheading {2. Effective estimates for pluricanonical embeddings for threefolds}

\sk

This section is motivated by the following

\proclaim{Question 2.1 }\ Let $Y$ be a smooth projective manifold of dimension $n$
which is of general type.
Does there exist an integer $m$, depending only on
$n$, such that the $m$-th pluricanonical map  $\Phi_{mK_Y}: Y \r \Phi_{mK_Y}(Y)$ is
birational onto its image?
\ep

It is well known (cf. \cite{B-P-V}) that for curves we can choose $m=3$, and
for surfaces we can choose $m=5$.
 Luo  conjectured in \cite{Luo1}, \cite{Luo2} that for the case of threefolds
 the answer to the question should also be affirmative.
In these two papers, he proves his conjecture in `almost all' possible
cases. Especially he  shows Theorem 0.3
 (cf. Theorem 5.1, Corollary 5.3 of \cite{Luo1}).

When the second named author gave a proof of Conjecture 0.2 
for threefolds (cf. \cite{Det}) he was not
aware of the papers \cite{Luo1} and \cite{Luo2} of Luo. 
So he independently
gave a proof of Theorem 0.3, using however the same basic idea  
(apparently both proofs were motivated by the paper \cite{Flet} of Fletcher).
Since the proof given in \cite{Det} seems to use the basic idea
in a shorter way and, moreover, easily gives effective bounds, we want to
include it here. More precisely we prove the following statement:

\proclaim{Theorem 2.2 }\ Let $C$ be a positive integer.
 Define $R=$lcm$(2,3,...$ $,26C-1)$ and
$m=$lcm$(4R+3,143C+5)$. Let $Y$ any  smooth projective threefolds  of
general type for
which $\chi (Y,\Cal {O}_Y) \leq C$ holds.
Then $\Phi_{mK_Y}: Y \r \Phi_{mK_Y}(Y)$ is birational onto its image.
\endproclaim

For the convenience of the reader and to fix further notations we 
recall some facts on which the proof is built.

We need the Plurigenus Formula due to Barlow, Fletcher and Reid
(cf. \cite{Flet}, \cite{Rei2}, see also \cite{Kol-Mor}, p.666 for the last part):

\proclaim{Theorem 2.3 }\  Let $Y$ be a projective threefold with only canonical singularities.
Then 
$$ \chi (Y, mK_Y) = \frac{1}{12}(2m-1)m(m-1)K_Y^3 - (2m-1) \chi (Y, 
\Cal {O}_Y)
+ \sum_Q l(Q,m)$$
with 
$$   l(Q,m) = \sum_{k=1}^{m-1} \frac{\overline{bk}(r- \overline{bk})}{2r}=
  \frac{r^2-1}{12}(m- \overline{m}) + 
 \sum_{k=1}^{\overline{m}-1} \frac{\overline{bk}(r- \overline{bk})}{2r}$$
Here the summation takes place over a basket of singularities $Q$ of type
$\frac{1}{r}(a,-a,1)$. $\overline{j}$ 
denotes the smallest nonnegative residue of $j$ modulo $r$, and $b$ is chosen such that 
$\overline{ab}=1$.

\noin Furthermore,
$$ \text{index}\  (Y) =  \text{lcm} \{r=r(Q): Q\in \text{basket} \} $$
\ep

 Hanamura (\cite{Han}) proves:

\proclaim{Theorem 2.4 }\
Let $Y$ be a smooth projective threefold of general type, which has a 
minimal or canonical model of index $r$. Then for any $m \geq m_0$ the $m$-th
pluricanonical map is birational onto its image, where
$$  m_0=4r+5  \text{ for }\  1 \leq r \leq 2$$
                 $$    m_0=4r+4 \text{ for }\ 3 \leq r \leq 5$$
                   $$  m_0=4r+3  \text{ for }\  r \geq 6 $$

\endproclaim

In the last step of the proof we use the following
theorem of Kollar (Corollary 4.8 in \cite{Kol2}):

\proclaim{Theorem 2.5 }\
Assume that for a smooth projective complex threefold $Y$ 
of general type we have $h^0(Y, lK_{Y}) \geq 2$. Then 
the $(11l+5)$-th pluricanonical map is birational onto its image.
\ep

 For estimating from below the terms $l(Q,m)$ in the Plurigenus 
Formula, we  need two Propositions due to Fletcher \cite{Flet}. In these Propositions
  $[s]$ denotes
the integral part of $s \in \BR$.

\proclaim{Proposition 2.6 }\
$$l( \frac{1}{r}(1,-1,1), m) = \frac{\overline{m} (\overline{m} -1)(
3r+1-2 \overline{m})}{12r} + \frac{r^2-1}{12} [\frac{m}{r}]$$
\ep

\proclaim{Proposition 2.7 }\ For $\alpha, \beta \in {\Bbb Z}$ with $0 \leq \beta \leq \alpha$ 
and for all $m \leq [(\alpha +1)/2]$, the following holds:
$$ l(\frac{1}{\alpha} (a,-a,1), m) \geq l(\frac{1}{\beta} (1,-1,1),m)$$
\ep

The basic idea of the proof is the following: We look at the
canonical model of the threefold $Y$, which exists by the famous
result of Mori \cite{Mor}, combined with results of Fujita \cite{Fuj},
 Benveniste \cite{Ben}
and Kawamata \cite{Kaw}. If the index of the canonical model is small, 
we can finish the proof by using Hanamura's Theorem. If the index  is
big,  we use the Plurigenus Formula due to Barlow, Fletcher and
Reid to show that for some $m$ we have $h^0(Y,mK_Y) \geq 2$,
and  finish the proof  by using Kollar's theorem.

\demo{ Proof of Theorem 2.2 }\
We first observe that by a theorem due to
 Elkik \cite{Elk} and Flenner
\cite{Flen} (cf. \cite{Rei2}, p.363),
canonical singularities are rational singularities.
  Hence, by the degeneration of the Leray spectral sequence
we have 
$$ \chi (Y, \Cal {O}_Y) = \chi (Y_c, \Cal {O}_{Y_c}).$$ 
If the index of $Y_c$ divides $R$, we  apply Hanamura's Theorem  and get that
$\Phi_{(4R+3)K_Y}$ embeds birationally.
 Hence, we may assume that the index
does not divide $R$. Then in the Plurigenus Formula
we necessarily have at least one singularity $\tilde{Q}$ in the basket of singularities
which is of the type $\frac{1}{r}(a,-a,1)$ with $r \geq 26C$.
Applying  a vanishing theorem for ample sheaves (cf. Theorem 4.1 in
\cite{Flet}),  the fact that $K^3_{Y_c} >0$ (since $K_{Y_c}$ is an ample ${\Bbb Q}$-divisor)
 and finally the  Propositions 2.6 and 2.7 of Fletcher, we obtain:  
$$ h^0(Y_c, (13C)K_{Y_c})$$ 
$$  =   \chi (Y_c, (13C)K_{Y_c})$$
$$  \geq   (1-26C)\chi (Y_c, \Cal {O}_{Y_c}) + \sum_{ Q \in \text{basket}} l(Q,13C)$$
$$  \geq   (1-26C)C + l(\tilde{Q},13C)$$
$$  \geq   (1-26C)C + l(\frac{1}{26C}(1,-1,1),13C)$$
$$  =     (1-26C)C + \frac{13C(13C-1)(78C+1-26C)}{312C}$$
$$  =     \frac{52C^2 - 15C - 1}{24}
 \geq  \frac{36}{24}  = 1.5 $$
\noin The last inequality is true since $C\geq 1$.
 Since $ h^0(Y_c, (13C)K_{Y_c})$ is an integer,
it has to be at least 2. 
 From the Definition of canonical singularities it easily follows
 (cf.
 e.g. \cite{Rei1}, p.277, \cite{Rei2}, p.355 or 
\cite{Flet}, p.225)
that  $h^0(Y,(13C)K_Y) \geq 2$. Now we can finish the proof 
by applying Theorem 2.5 due to Kollar.\qed \edm

\noindent Despite the fact that our $m=m(C)$ is explicit, it is so huge that it
is only of theoretical interest. For example for $C=1$ 
 one can choose $m=269$ (\cite{Flet}), but for $C=1$ our $m$ is
already for of the size $10^{13}$. Moreover, for all examples
of threefolds of general type which are known so far,
any $m \geq 7$ works. So we guess there should exist a bound which is 
 independent of the size of the holomorphic Euler
characteristic.
   
\bigskip

\subheading{3. Iitaka-Severi's Conjecture for threefolds}

\sk

 The claim of this section is the following 
\proclaim{Theorem 3.1 }\
Let $X$ be a fixed complex variety. Then  the set of targets
 $\Cal {F}(X)/ \sim$ 
 with   dim$_{\BC}Y \leq 3$ is a finite set.
\ep
 
 By Proposition 6.5 of Maehara \cite{Mae2} it is sufficient to show 
the following: 
There exists a natural number $m$, only depending on $X$, 
such that $\Cal {F}(X) \s \Cal {F}_m(X)$ for varieties $Y$
with dim$_{\BC}Y \leq 3$. 
Since we  prove finiteness only up to birational equivalence,
we may assume, without loss of generality, that $X$ and all $Y$ in Theorem 3.1
 are nonsingular projective varieties. This is by virtue
of Hironaka's resolution theorem \cite{Hir}, cf. also \cite{Uen}, p.73. 
Hence, using Theorem 2.2 or Theorem 0.3 of Luo we get Theorem 3.1 
as a consequence of  the following: 

\sk
\proclaim{Proposition 3.2 }\ 
Let $X$ be a fixed smooth projective variety and $f:X \r Y$ a dominant
 rational map
to another smooth projective variety $Y$ with dim$_{\BC}Y =n$.
Then we have $$\chi (Y, \Cal {O}_Y) \leq \sum_{\{i|2i \leq n\}}
 h^i(X,\Cal {O}_X) $$
\ep \sk

\demo{ Proof of Proposition 3.3 }\
  First we obtain, 
  by Hodge theory on compact K\"ahler manifolds 
  (cf. \cite{Gri-Har}, or \cite{Iit}, p.199) 
  $$ h^i(Y, \Cal {O}_Y) = h^0(Y, \Omega_Y^i),$$
where $i=1,...n$. The same kind of equalities hold for $X$.
  Now by \cite{Iit}, Theorem 5.3, we obtain that
  $$ h^0(Y, \Omega_Y^i) \leq h^0(X, \Omega_X^i),$$
where again $i=1,...,n$. Hence, we can conclude: 
$$\chi (Y, \Cal {O}_Y) \leq
\sum_{\{i|2i \leq n\}}
 h^i(Y,\Cal {O}_Y) = 
\sum_{\{i|2i \leq n\}}
 h^0(Y,\Omega_Y^i) \leq$$ $$ \leq
\sum_{\{i|2i \leq n\}}
 h^0(X,\Omega_X^i) =
\sum_{\{i|2i \leq n\}}
 h^i(X,\Cal {O}_X) 
 $$ \qed \edm

\bigskip

\subheading {4.  On the number of Targets}

\sk

Let $X$ be a smooth threefold of general type, and define $r=r_{X_c}$, $k=K_{X_c}^3$.

\sk
\proclaim {Theorem 4.1 }\ 
There exists a universal constant $C(r,k)$, depending only on $r$ and $k$,
such that
$$\#( \Cal {F}(X) / \sim ) \leq C(r,k) $$
 \ep

\sk

\sk

\proclaim {Theorem 4.2}
 There exists a universal constant $C'(r,k)$, depending only on $r$ and $k$,
such that
 if $Y$ is a smooth threefold and $R(X,Y) \neq \emptyset$,
 then
$$r_{Y_c} \leq C'(r,k) $$
\endproclaim
\sk

The rest of this section deals with the proof of Theorem 4.1 and Theorem 4.2,
which we prove simultaneously. 
We fix positive integers $r$ and $k$.
Denote by $\BX(r,k)$ the set of  threefolds $X_c$ with only canonical
singularities and ample canonical sheaves $K_{X_c}$, which satisfy 
 $r_{X_c}=r$, $K_{X_c}^3=k$. Let $X$ be a smooth threefold such that $X_c \in \BX(r,k)$. 

   a) In this part of the proof we only consider targets 
$((Y,f)/ \sim ) \in (\Cal {F}(X) / \sim) $ with dim$_{\BC}Y =3$.

   Due to Theorem 2.4 of Hanamura,  the map  
 $$\Phi _{9rK_X}:X \r \BP^N $$
is birational onto its image,
 where, by ([Mat-Mum])
  $$N=\dim H^0(X, 9rK_X)-1 \leq 9^3 r^3 k+3.$$ 
 Moreover, by [Ban-Mar],  Lemma 1 (cf. also Proposition 2, part 2), the degree $ d_X$ 
of the image $X'=
\Phi_{9rK_X} (X)$ has the bound 
$$d_X \leq 9^3r^3 k.$$ 

Let $Y$ be a smooth threefolds of general type with $R(X,Y)
\not= \emptyset$. 

\sk
\proclaim{Proposition 4.3 }\
There exists a universal constant $C_1(r,k)$, depending only on $r$ and $k$, 
such that  we have
$$\chi (Y, \Cal {O}_Y) \leq C_1(r,k)$$
\ep
\sk

\demo{ Proof of Proposition 4.3 }\
By Proposition 3.2 we have
$$\chi (Y,\Cal {O}_Y) \leq h^2(X, \Cal {O}_X) +1.$$
In the Hilbert polynomials $\chi (X_c, mrK_{X_c})$  the expressions
$l(Q, rm)$, cf. 
Theorem 2.3, are linear in $m$, and so the two highest
 coefficients of the polynomial $\chi (X_c, mrK_{X_c})$ in the variable
 $m$ only depend on $r^3k$.
But then by Theorem 2.1.3 of Kollar \cite{Kol1}, the family of  the 
$(X_c, rK_{X_c})$, where $X_c \in \BX(r,k)$,  is a bounded family. 
That means there exists a morphism
$\pi : \Cal {X} \r S$ between (not necessarily complete) varieties $\Cal {X}$ and $S$ and a $\pi$-ample
Cartier divisor $D$ on $\Cal {X}$ such that every $(X_c, rK_{X_c})$ is
isomorphic to  $(\pi^{-1}(s), D|_{\pi^{-1}(s)})$ for some $s \in S$.
So it is sufficient  to prove that there exists a constant $C_0$ which satisfies:
For all $s \in S$ and for some desingularization $X(s)$ of $\pi^{-1}(s)$
we have $h^2(X(s),\Cal {O}_{X(s)}) \leq C_0$\footnote{Remark that by an easy argument
like in the proof of Proposition 3.2 any two such desingularizations have the same
$h^2(X(s),\Cal {O}_{X(s)})$.}.

This is  shown by using first generic uniform desingularization of
the family $ \pi : \Cal {X} \r S$ (cf. \cite{Hir}, \cite{Bin-Fle}), and
 afterwards a semi-continuity theorem (cf. \cite{Gro}, \cite{Gra}):
By applying generic uniform desingularization and induction on the dimension
there exist finitely many subvarieties $S_i, i=1,...,l$, which cover $S$,
and morphisms $\Psi_i : \Cal {Y}_i \r S_i$ between varieties $\Cal {Y}_i$ and
 $S_i$ which desingularize $\Cal {X}_i := \pi^{-1}(S_i)$ fiberwise, i.e., there exist
morphisms $\Phi_i :\Cal {Y}_i \r \Cal {X}_i$ over $S_i$ such that for any 
$s \in S_i$ the map $\Phi_i :\Psi_i^{-1}(s) \r \pi^{-1}(s)$ 
is a desingularization.

Using semi-continuity for the families $\Phi_i : \Cal {Y}_i \r S_i$,
we  obtain finitely many subvarieties $S_{ij}, j=1,...,l_i$ of $S_i$,
which cover $S_i$, and have the following property: 
If we denote  $\Cal {Y}_{ij}:= \Phi_i^{-1}(S_{ij})$ and $\Phi_{ij}:= 
\Phi_i|_{\Cal {Y}_{ij}}$,
we get that for  the families $\Phi_{ij}: \Cal{Y}_{ij} \r S_{ij}$ the number 
$C_{ij}:= h^2( \Phi_{ij}^{-1}(s),\Cal {O}_{\Phi_{ij}^{-1}(s)})$ is constant for 
 $s \in S_{ij}$.
Hence,  $C_0:= \max_{i=1,...,l;j=1,...,l_i} C_{ij}$ has the desired property. 
\qed \edm
\sk
\noindent {\bf Remark.}  Proposition 4.3 can also be proved
as follows. By a result of Milnor ([Mil])  the Betti numbers 
of the variety $X'=\Phi_{9rK_X}(X)$ have estimates depending on its
 degree $d_X \leq 9^3r^3 k$
  only. From the standard exact cohomology
sequences and dualities it easily follows that $h^{2,0}(X)$ may be 
estimated by Betti numbers of $X'$. 
\sk

Using Proposition 4.3 and Theorem 2.2 we can choose an integer $p=p(r,k)$,  such that
$p$ is divisible by $r$, $p \geq 9r$ and
 $$\Phi _{pK_Y}:Y \r \BP^M $$
is birational onto its image,
 where, by ([Mat-Mum])
  $$M=\dim H^0(Y, pK_Y)-1 \leq p^3k+3.$$ 

\proclaim {Lemma 4.4}

1) The degree of $Y' :=\Phi_{pK_Y} (Y)\subset \BP^M$ 
is smaller than deg$X' \leq p^3k$.

2) For any map $f\in R(X',Y')$ the degree $d_f$
of  its graph $\Gamma _f \subset \BP^N \times \BP^M$ is 
not greater than $8p^3k$.
\ep
Lemma 4.4 is a particular case  of part 2 and 3 of
Proposition 2 of \cite {Ban-Mar} for  $n=3$, applied to  
 the threefolds $X_c, Y_c$ and 
linear systems $|pK_{X_c}|,|pK_{Y_c}|$. 
We have to note only that  Proposition 2 and Lemma 1
in \cite {Ban-Mar} 
is stated for Cartier divisors. But
only the fact that they are ${\Bbb Q}$-Cartier is used in their proofs. \qed

\sk

 By Proposition 1 of the same paper ([Ban-Mar]),
there exist  algebraic families $(\CX,p_X,T)$, $(\CZ,p_Z,V)$, $(\CY,p_Y,U)$
with constructive bases  and projections
$\pi_U:V \to U, \pi_T: V \to T$,    with the following properties:

1) For any $X_c \in \BX(r,k)$, there is a point $t \in T$, such that
$X_c$ is birational to $X=p_X^{-1}(t)$, and all points $t \in T$ have this property.

2) For any $Y$ with $R(X,Y) \neq\emptyset$ for some $X_c\in \BX(r,k)$,
there is a point $u\in U$, such that $Y$ is birational
to $Y =p_Y^{-1}(u)$, and all points $u \in U$ have this property.

3) For any dominant rational map $f:p_X^{-1}(t)=X \to Y=p_Y^{-1}(u)$,
 there is a point $v\in V$, such that 
$\pi_U(v)=u $, $\pi_T(v)=t$, $p_V^{-1}(v)$ is a graph of the map $f$, and all points
$v \in V$ have this property.

Let $\tilde{V} =\{ (t,v) \in T \times V| \pi_T (v)=t \}$, and denote by $p_T$ resp.
$p_V$ the projections to the first resp. to the second factor. Through the composed map
$\pi_U \circ p_V : \tilde{V} \r U$ the variety $\tilde{V}$ is also a variety over $U$.
Let $\tilde{\CY}= \tilde{V} \times_U \CY$ be obtained by base change, and denote the
projection to the first factor by $p_{\tilde{V}}: \tilde{\CY} \r \tilde{V}$.
Then we have
$$ \alignat5 
& \tilde{\CY}   \quad && \overset   {p_{\tilde{V}}} \to \longrightarrow \quad&& \tilde{V} \quad&&
\overset   {p_{T}} \to \longrightarrow \quad && T.\\
\endalignat$$
In this diagram, for every $ t \in T$, the set $p_T^{-1}(t)$ can be considered as the set
of graphs of dominant rational maps $f: X \r Y$, where $X=p_X^{-1}(t)$, and $p_{\tilde{V}}:
\tilde{\CY} \r \tilde{V}$ is the universal family of threefolds $Y$ over the graphs
of $f:X \r Y$.

By applying the process of local uniform desingularization, described in Proposition 4.3,
to the family  $p_{\tilde{V}}: \tilde{\CY} \r \tilde{V}$,
 we obtain a finite number of smooth families
 $(p_{\tilde{V}})_i: (\tilde{\CY})_i \r (\tilde{V})_i$, $i=1,...,l$, 
 the bases $(\tilde{V})_i$ 
of which are connected and  cover $\tilde{V}$, and the fibers of which 
are desingularizations of 
the fibers of $p_{\tilde{V}}: \tilde{\CY} \r \tilde{V}$.
 For any $i$ the map
 $(p_{\tilde{V}})_i: (\tilde{\CY})_i \r (\tilde{V})_i$ 
is a smooth family of projective threefolds of general type over a connected
 base 
$(\tilde{V})_i$. By a theorem of J.Kollar and Sh.Mori (\cite{Kol-Mor},
 Theorem 12.7.6.2)
there is an algebraic map $\phi_i$ from $(\tilde{V})_i$ to the birational
equivalence classes of the fibers of 
$(p_{\tilde{V}})_i: (\tilde{\CY})_i \r (\tilde{V})_i$. Moreover, all these
fibers have the same Hilbert function.

\noindent From this fact two conclusions can be derived:

{\bf 1.} Since the index of a canonical threefold can be bounded in terms
of the Hilbert function (\cite{Kol-Mor}, p.666), the indices of the canonical
models of the fibers of the family $(p_{\tilde{V}})_i: 
(\tilde{\CY})_i \r (\tilde{V})_i$
vary in a finite set of natural numbers, only.

 {\bf 2.} Let $(p_T)_i:=p_T|_{(\tilde{V})_i}$ and  $n_i(t)$ be the number of 
irreducible components of $(p_T)_i^{-1}(t)$
(it may be zero). Define, for  $X = p_X^{-1}(t)$, $\Cal {G} (X) = \{ Y|(Y,f) \in \Cal {F}(X) \}$, 
and let $\sim$
denote birational equivalence on $\Cal {G} (X)$. 
Since $(\# \Cal {G}(X) / \sim) < \infty$, it follows that
the restriction  $\phi_i$ to  $(p_T)_i^{-1}(t)$ has to be constant
on the connected components of  $(p_T)_i^{-1}(t)$.
Then
$$ \#( \Cal {G}(X) / \sim) \leq \sum_{i=1}^l n_i(t). $$

Since from the beginning  the constructions of all the families
 were algebraic
and defined only by the constants $r$ and $k$, we have proved Theorem 4.2, and 
also
the following
\sk
\proclaim {Lemma 4.5}
There exists a universal constant $C_2(r,k)$, depending only on $r$ and $k$,
such that we have
$$(\# \Cal {G}(X)/ \sim ) \leq C_2(r,k).$$
\ep
Next, we look at the map 
$$ (\pi_T, \pi_U) : V \r T \times U.$$
It is algebraic and any point of the fiber over $(t,u)\in T \times U$
 defines a map from
$R(\pi_T^{-1}(t),\pi_U^{-1}(u))$. The last set is finite for all $(t,u)
 $, so we get:
\proclaim {Lemma 4.6}
There exists a universal constant $C_3(r,k)$, depending only on $r$ and $k$,
such that 
$$\# R(X,Y) \leq C_3(r,k).$$
\ep

 From Lemma 4.5 and Lemma 4.6 the statement of Theorem 4.1 for
3-dimensional targets  is immediate.
The desired bound may be chosen as
$C_2(r,k)C_3(r,k) $.\hfill \qed

b) 
Now we consider targets 
$((Y,f)/ \sim) \in (\Cal {F}(X) / \sim)$ with dim$_{\BC}Y \leq 2$.
For these targets we know that the indices
 of the $Y_c$ are 1 or 2. So we can repeat the same argument as above, omitting however
Proposition 4.3. The only change which has to be done is 
 replacing the moduli spaces due to [Kol-Mor]
by the respective moduli spaces for surfaces or curves. So we get Theorem 4.1, and
in particular   
Lemma 4.5 and Lemma 4.6, also for these kinds of targets. 
\hfill \qed
\sk
\noindent {\bf Remark.}
According to \cite{Ban-Mar}, there exists a universal function $\sigma$ in two
variables, such that $\#R(X,Y) \leq \sigma (r \cdot r_{Y_c},k)$. This fact,
together with Theorem 4.2, yields an alternative proof of Lemma 4.6.

\bigskip

\subheading{5. A Conjecture of Kobayashi for threefolds of general type}

\sk

 In this section we prove 
 \sk
 
\proclaim {Theorem 5.1 }\ For any complex variety  $X$
there is a number $c(X)$ such that 
         $$\#R(X,Y) \leq c(X)$$
 for any  complex variety $Y$ of general type with dim$_{\BC}Y \leq 3$.
If  $X$ is a threefold of general type, then $c(X)$
can be expressed only in terms of $r_{X_c}$ and  $K_{X_c}^3$. 
\endproclaim
\sk

\demo{Proof of Theorem 5.1}\ 
Like in section 3 we may assume that $X$ and $Y$ are smooth
projective varieties. By \cite{Kob-Och},  $\#R(X,Y)$ is finite for
 every fixed $Y$. By Theorem 3.1, we know that
  for given $X$ there exist only
 finitely many such $Y$, up to birational 
equivalence. Since birational equivalence does not
effect the number $\#R(X,Y)$, the first statement follows.

Let  $X$ now be a projective threefold of general type. 
 Then the second statement is just Lemma 4.6. \qed \enddemo

\sk
\proclaim{Remark 5.2}
The estimate which is given in Theorem 5.1 is not effective.
\endproclaim

\bk

\noindent Tatiana M. Bandman, Dept. of Mathematics and Computer Sciences, 
 Bar-Ilan University, Ramat-Gan, Israel, e-mail: bandman\@macs.cs.biu.ac.il

\noindent Gerd Dethloff, Dept. of Mathematics, Osaka University, Toyonaka, Osaka 560, 
Japan, e-mail:dethloff\@math.sci.osaka-u.ac.jp

\end